\documentclass[12pt, a4paper,reqno]{amsart}

\usepackage{amsmath}
\usepackage{amsfonts}
\usepackage{amssymb}
\usepackage{amsthm}
\usepackage[latin1]{inputenc}

\numberwithin{equation}{section}

\theoremstyle{plain}
\newtheorem{theorem}[equation]{Theorem}
\newtheorem{corollary}[equation]{Corollary}
\newtheorem{lemma}[equation]{Lemma}

\theoremstyle{definition}

\newtheorem{remark}[equation]{Remark}

\numberwithin{equation}{section}

\def\mvint_#1{\mathchoice%
          {\mathop{\kern 0.2em\vrule width 0.6em height 0.69678ex depth -0.58065ex
                  \kern -0.8em \intop}\nolimits_{\kern -0.4em#1}}%
          {\mathop{\kern 0.1em\vrule width 0.5em height 0.69678ex depth -0.60387ex
                  \kern -0.6em \intop}\nolimits_{#1}}%
          {\mathop{\kern 0.1em\vrule width 0.5em height 0.69678ex depth -0.60387ex
                  \kern -0.6em \intop}\nolimits_{#1}}%
          {\mathop{\kern 0.1em\vrule width 0.5em height 0.69678ex depth -0.60387ex
                  \kern -0.6em \intop}\nolimits_{#1}}}

\newcommand{\R}{\mathbb{R}}

\newcommand\dist{\operatorname{dist}}

\newcommand\supp{\operatorname{supp}}
\newcommand\vlc{\operatorname{VLC}}
\begin{document}

\title[How to recognize polynomials]{How to recognize polynomials in higher order Sobolev spaces}
\author[Bojarski, Ihnatsyeva \and Kinnunen]{Bogdan Bojarski, Lizaveta Ihnatsyeva \and Juha Kinnunen }
\thanks{The research is supported by the Academy of Finland, the first author is also partially supported by Polish Ministry of Science grant no N N201 397837 (years 2009-2012).}

%
%\keywords{XXXXX}
%%%%%%%%%%%%%%%%%%%%%%%%%%%%%%%%%%%%%%
%%%%%%%%%%%%%%%%%%%%%%%%%%%%%%%%%%%%%%
\begin{abstract}
This paper extends characterizations of Sobolev spaces by Bourgain, Br\'{e}zis, and Mironescu to the higher order case.
As a byproduct, we obtain an integral condition for the Taylor remainder term, which implies that the function is a polynomial.
Similar questions are also considered in the context of Whitney jets.
\end{abstract}
\maketitle

\section{Introduction}
%\footnotetext[1]{The research is supported by Academy of Finland, the first author is also partially supported by Polish Ministry of Science grant no N N201 397837 (years %2009-2012).}

\medskip
In this paper we study a new characterization of the higher order Sobolev spaces $W^{m,p}(\Omega)$ which is based on J. Bourgain, H. Br\'{e}zis, and P. Mironescu's approach \cite{BourgainBrezisMironescu} (see also \cite{Brezis}).
They showed that a function $f\in L^p(\Omega)$ belongs to the first order Sobolev space $W^{1,p}(\Omega)$, ${1<p<\infty}$, on a smooth bounded domain $\Omega\subset\mathbb{R}^n$ if and only if
\begin{equation}\label{BrezisCondition}
\liminf\limits_{\varepsilon\to 0}\int_{\Omega}\int_{\Omega}\frac{|f(x)-f(y)|^p}{|x-y|^{p}}\rho_\varepsilon(|x-y|)\,dx\,dy<\infty,
\end{equation}
where $\rho_\varepsilon$, with $\varepsilon>0$, are radial mollifiers.
Moreover,
\[
\liminf\limits_{\varepsilon\to 0}\int_{\Omega}\int_{\Omega}\frac{|f(x)-f(y)|^p}{|x-y|^{p}}\rho_\varepsilon(|x-y|)\,dx\,dy
=c\int_\Omega|\nabla f|^p\,dx,
\]
where the constant $c$ depends only on $p$ and $n$.
For $p=1$ this gives a characterization of the space of bounded variation $BV(\Omega)$.
See also \cite{BourgainNguyen}, \cite{Davila}, \cite{Nguyen06}, \cite{Nguyen08} and \cite{Ponce} for related results.

We extend the results of \cite{BourgainBrezisMironescu} and \cite{Brezis}
to the higher order case.
To characterize the Sobolev spaces $W^{m,p}(\Omega)$, ${1<p<\infty}$, we use the condition
\begin{equation}\label{eqCondition}
\liminf\limits_{\varepsilon\to 0}\int_{\Omega}\int_{\Omega}\frac{|R^{m-1}f(x,y)|^p}{|x-y|^{mp}}
\rho_\varepsilon(|x-y|)\,dx\,dy<\infty,
\end{equation}
where $R^{m-1}f$ is a Taylor $(m-1)$-remainder of $f$, generalizing \eqref{BrezisCondition}.
For $p=1$ the condition \eqref{eqCondition} describes the space $BV^m(\Omega)$ of integrable functions whose
weak derivatives of order $m$ are signed Radon measures with finite total variation.

Condition \eqref{eqCondition} is a priori weaker than the pointwise condition
\begin{equation}\label{PointwiseCharacteristic}
|R^{m-1}f(x,y)|\le |x-y|^m(a_f(x)+a_f(y)),\,\,\,a_f\in L^p(\R^n),
\end{equation}
characterizing Sobolev class $W^{m,p}(\R^n)$ as in \cite{BojarskiProcSteklovInst} (see also (\ref{VLCcondition}) below).

Another variant of extension of the results of \cite{BourgainBrezisMironescu} and \cite{Brezis} to the higher order case
has been introduced in \cite{Borghol}, where the characterization of $W^{m,p}(\Omega)$, $1<p<\infty$, ($BV^m(\Omega)$ for $p=1$) is formulated in terms of the $m$-th differences.
According to this result, a function $f\in L^p(\Omega)$ belongs to $W^{m,p}(\Omega)$ if
%the formal Taylor
%$(m-1)$-remainder satisfies
%\[
%|R^{m-1}f(x,y)|\le|x-y|[a(x)+a(y)]
%\]
%outside a set of measure zero for some function $a\in L^p(\Omega)$.
\begin{equation}\label{CharactViaDifferences}
\liminf\limits_{\varepsilon\to 0}\int_{\Omega}\int_{\Omega}\bigg|\sum\limits_{j=0}^{m}(-1)^j\binom{m}{j}f\bigg(\frac{(m-j)}{m}x+\frac{j}{m}y\bigg)\bigg|^p
\frac{\rho_\varepsilon(|x-y|)}{|x-y|^{mp}}\,dx\,dy<\infty.
\end{equation}

For smooth functions $f\in C^{m+1}$ the equivalence of the integrands in \eqref{eqCondition} and \eqref{CharactViaDifferences} modulo $O(|x-y|^{m+1})$ is well known. The results of the present paper and of \cite{Borghol} essentially show that both integrands are equivalent in their averaged asymptotic behaviour, for $\varepsilon\to 0$, in the $\varepsilon$-neighbourhood of the diagonal $\Delta=\{x=y\}$ in the Cartesian product $\Omega\times\Omega$.

In close connection with these  characterizations,
H. Br\'{e}zis \cite{Brezis} considered  conditions under which a measurable function
$f$ defined on a connected open set $\Omega$  is a constant. See also \cite{Pietrushka}. In particular,
he showed that if
\[
\int_\Omega\int_\Omega\frac{|f(x)-f(y)|}{|x-y|^{n+1}}\,dx\,dy<\infty,
\]
then $f$ is a constant function.
We extend this result to the higher order case and show that the condition
\[
\int_\Omega\int_\Omega\frac{|R^{m-1}f(x,y)|}{|x-y|^{n+m}}\,dx\,dy<\infty
\]
implies that the function $f$, with locally integrable weak derivatives up to order $m-1$, is a polynomial of degree at most $m-1$.

The condition (\ref{eqCondition}) applies to Whitney jets as well. Recall that H. Whitney in \cite{Whitney} gave a method to define differentiable functions on closed subsets of $\R^n$.
His approach can be adopted to different kind of smoothness conditions.
In particular, for an $(m-1)$-jet $F$ on a subset $\Omega\subset\R^n$, defined as a collection of functions $\{f_\alpha:\,|\alpha|\le m-1\}$, $f_\alpha\in L^p(\Omega)$, we can study the meaning of (\ref{eqCondition}) using the formal Taylor remainder of order $m-1$ of the jet $F$.
The formalism of Taylor-Whitney jets identifies in a natural way virtual derivatives, jets, with Sobolev derivatives.

In Section \ref{sWhitneyjets} we show that a jet $F$ on an open set $\Omega$ satisfying (\ref{eqCondition}) for some special case of mollifiers is locally a jet of $W^{m,p}$-function. The proof uses an approximation procedure from \cite{BojarskiProcSteklovInst} where the Sobolev spaces are described in terms of pointwise inequalities \eqref{PointwiseCharacteristic}.

Notice also that a certain version of condition (\ref{eqCondition}) for the jet $F$ on Ahlfors $s$-regular subsets $S\subset\R^n$, $n-1<s\le n$, charaterizes the Lipschitz spaces ${\rm Lip }(m,p,\infty,S)$ studied in \cite{JonssonWallin}. If $S$ supports the $q$-Poncar\'{e} inequality, $1\le q<\infty$, then the first order space ${\rm Lip }(1,p,\infty,S)$ coincides with the Haj\l asz-Sobolev space $W^{1,p}(S)$ for $p>q$ (see e.g. \cite{YangLin}).

\section{Characterization of Sobolev spaces}

Our notation is standard.
For a multi-index $\alpha=(\alpha_1,\dots,\alpha_n)$, $\alpha_i\geq 0$, and a point
$x=(x_1,\dots,x_n)\in\R^n$, we denote by
\[
x^\alpha=x_1^{\alpha_1}x_2^{\alpha_2}\cdots x_n^{\alpha_n}
\]
the monomial of degree
\[
|\alpha|=\sum\limits_{i=1}^n \alpha_i.
\]
In the same way
\[
D^{\alpha}f=\dfrac{\partial^{|\alpha|}f}{\partial x_1^{\alpha_1}\dots
\partial x_n^{\alpha_n}}
\]
is a (weak) partial derivative of order $|\alpha|$. We also use the convention that $D^{0}f=f$.
Moreover, let $\nabla^mf$ be a vector with the components $D^\alpha f$, $|\alpha|=m$.

Let $\Omega$ be an open set in $\R^n$, $1\le p<\infty$, $m$ a positive integer.
The Sobolev space $W^{m,p}(\Omega)$ consists of all functions $u\in L^p(\Omega)$ such that
for all multi-index $\alpha$ with $|\alpha|\le m$ the weak derivative $D^\alpha u$ exists and
belongs to $L^p(\Omega)$. We use the convention $W^{0,p}(\Omega)=L^p(\Omega)$.
The Sobolev space $W^{m,p}(\Omega)$ is equipped with the norm
\[
\|u\|_{W^{m,p}(\Omega)}=\Big(\sum_{|\alpha|\le m}\int_\Omega|D^\alpha u|^p\,dx\Big)^{1/p}.
\]
For the properties of Sobolev functions, see \cite{Mazya}.

We write
\[
T_y^{m}f(x)=\sum_{|\alpha|\le m}D^{\alpha}f(y)\frac{(x-y)^\alpha}{\alpha!}
\]
and
\[
R^{m}f(x,y)=f(x)-T_y^{m}f(x)
\]
for the Taylor polynomial of order $m$ and the Taylor remainder of order $m$, respectively.
We use the same notation for formal Taylor polynomials and remainders if we have only weak derivatives.
Let $\Omega$ be an open set in $\R^n$, $1\le p<\infty$, $m$ a positive integer.

We define a family of functions $\rho_\varepsilon\in L^{1}_{\rm loc}(0,\infty)$, $\varepsilon>0$,
such that $\rho_\varepsilon\geq 0$,
\[
\int_0^\infty\rho_\varepsilon(r)r^{n-1}dr=1
\]
and
\[
\lim\limits_{\varepsilon\to 0}\int_\delta^\infty\rho_\varepsilon(r)r^{n-1}dr=0\quad\text{for every}\quad\delta>0.
\]
These properties are rather standard in the construction of radial
mollifiers related to approximations of unity.

First we prove a useful result for smooth functions.

\begin{lemma}\label{LemmaEqualityForSmoothFunction}
Let $\Omega$ be an open set in $\R^n$, $1\le p<\infty$, $m$ a positive integer and $g\in C^{m+1}_0(\R^n)$.
Then
\begin{equation}\label{eqLemma}
\begin{split}
\lim_{\varepsilon\to 0}\int_{\Omega}\int_{\Omega}&\frac{|R^{m-1}g(x,y)|^p}{|x-y|^{mp}}\rho_\varepsilon(|x-y|)\,dx\,dy
\\
&=\int_{\Omega}\int_{\partial B(0,1)}\Big|\sum_{|\alpha|=m}\frac{D^{\alpha}g(x)}{\alpha! }e^\alpha\Big|^p\,de\,dx.
\end{split}
\end{equation}
\end{lemma}

\begin{proof}
Since $g\in C^{m+1}_0(\R^n)$  by Taylor's formula, we have
$$
|R^{m-1}g(x+h,x)|\le \Big|\sum_{|\alpha|=m}\frac{D^{\alpha}g(x)}{\alpha!}h^\alpha\Big|+c'|h|^{m+1}
$$
from which we conclude that
$$
|R^{m-1}g(x+h,x)|^p\le (1+\theta)
\Big|\sum_{|\alpha|=m}\frac{D^{\alpha}g(x)}{\alpha!}h^\alpha\Big|^p+c'_{\theta}|h|^{(m+1)p}
$$
for every $\theta>0$, $x\in\R^n$ and $h\in \R^n$.

We multiply the last inequality by $\rho_\varepsilon(|h|)/|h|^{mp}$ and integrate over the set
\[
\begin{split}
S=\{(x,h)&\in (\supp g\cap \Omega)\times\R^n:\, x+h\in \Omega\}
\\
&\cup\{(x,h)\in(\Omega\setminus\supp g)\times\R^n:\, x+h\in (\supp g\cap\Omega)\}.
\end{split}
\]

We have
\begin{equation}\label{InequalityLemma}
\begin{split}
\iint_S&\frac{|R^{m-1}g(x+h,x)|^p}{|h|^{mp}}\rho_\varepsilon(|h|)\,dh\,dx
\\
&\le(1+\theta)\int_{\Omega}\int_{\mathbb{R}^n}\frac{\rho_\varepsilon(|h|)}{|h|^{mp}}
\Big|\sum_{|\alpha|=m}\frac{D^{\alpha}g(x)}{\alpha!}h^\alpha\Big|^p\,dh\,dx
\\
&\qquad+2c'_\theta|{\rm supp}\,g|\int_{\mathbb{R}^n}|h|^{p}\rho_\varepsilon(|h|)\,dh.
\end{split}
\end{equation}

By the properties of the mollifiers $\rho_\varepsilon$, it follows that
\begin{equation}\label{Estimate}
\lim_{\varepsilon\to 0}\int_{\mathbb{R}^n}|h|^{p}\rho_\varepsilon(|h|)\,dh=0.
\end{equation}

Note also, that
\[
\begin{split}
\int_{\mathbb{R}^n}&\frac{\rho_\varepsilon(|h|)}{|h|^{mp}}
\Big|\sum_{|\alpha|=m}\frac{D^{\alpha}g(x)}{\alpha!}h^\alpha\Big|^p\,dh
\\
&=\int_{0}^\infty r^{n-1}\rho_\varepsilon(r)\,dr\int\limits_{\partial B(0,1)}
\Big|\sum_{|\alpha|=m}\frac{D^{\alpha}g(x)}{\alpha! }e^\alpha\Big|^p\,de.
\end{split}
\]

Thus, passing to the limit in  (\ref{InequalityLemma}), first with $\varepsilon\to 0$, then with $\theta\to 0$,
and changing variables in the integral on the left hand side, we arrive at
\begin{equation}\label{eqUpperEstimateLemma}
\begin{split}
\limsup_{\varepsilon\to 0}\int_\Omega\int_{\Omega}&\frac{|R^{m-1}g(x,y)|^p}{|x-y|^{mp}}\rho_\varepsilon(|x-y|)\,dx\,dy
\\
&\le\int_{\Omega}\int_{\partial B(0,1)}\Big|\sum_{|\alpha|=m}\frac{D^{\alpha}g(x)}{\alpha! }e^\alpha\Big|^p\,de\,dx.
\end{split}
\end{equation}
This concludes the first part of the proof of \eqref{eqLemma}.

Then we show the inequality to the other direction.
If  $K$ is a compact subset of $\Omega$, then for any $x\in K$ and
$|h|\le d$, where
\[
d=\min\{1,{\rm dist}\,(K,\partial\Omega)/2\},
\]
we have
$$
\Big|R^{m-1}g(x+h,x)-\sum_{|\alpha|=m}\frac{D^{\alpha}g(x)}{\alpha!}h^\alpha\Big|\le c_K|h|^{m+1}.
$$
Hence
$$
\Big|\sum_{|\alpha|=m}\frac{D^{\alpha}g(x)}{\alpha!}h^\alpha\Big|^p\le (1+\theta)|R^{m-1}g(x+h,x)|^p+c_{\theta,K}|h|^{(m+1)p}
$$
for every $\theta>0$ and consequently
\begin{equation}\label{BasicInequalityForProof}
\begin{split}
\int_K\int_{B(0,d)}&\frac{\rho_\varepsilon(|h|)}{|h|^{mp}}
\Big|\sum_{|\alpha|=m}\frac{D^{\alpha}g(x)}{\alpha!}h^\alpha\Big|^p\,dh\,dx
\\
&\le (1+\theta)\int_K\int_{B(0,d)}\frac{|R^{m-1}g(x+h,x)|^p}{|h|^{mp}}\rho_\varepsilon(|h|)\,dh\,dx
\\
&\qquad+c_{\theta,K}|K|\int_{B(0,d)}|h|^{p}\rho_\varepsilon(|h|)\,dh.
\end{split}
\end{equation}

Passing to the limit as $\varepsilon\to 0$ in (\ref{BasicInequalityForProof}),
and taking into account (\ref{Estimate}), we have
\[
\begin{split}
\int_K&\int_{\partial B(0,1)}
\Big|\sum_{|\alpha|=m}\frac{D^{\alpha}g(x)}{\alpha!}e^\alpha\Big|^p\,de\,dx
\\
&\le (1+\theta)\liminf_{\varepsilon\to 0}\int_{K}\int_{B(0,d)}\frac{|R^{m-1}g(x+h,x)|^p}{|h|^{mp}}\rho_\varepsilon(|h|)\,dh\,dx
\\
&\le(1+\theta)\liminf_{\varepsilon\to 0}\int_\Omega\int_{\Omega}\frac{|R^{m-1}g(x,y)|^p}{|x-y|^{mp}}\rho_\varepsilon(|x-y|)\,dx\,dy.
\end{split}
\]

Since the last estimate holds for every $\theta>0$ and every compact set $K\subset\Omega$, we have
\begin{equation}\label{eqLowerEstimateLemma}
\begin{split}
\int_{\Omega}\int_{\partial B(0,1)}&
\Big|\sum_{|\alpha|=m}\frac{D^{\alpha}g(x)}{\alpha!}e^\alpha\Big|^p\,de\,dx
\\
&\le\liminf_{\varepsilon\to 0}\int_\Omega\int_{\Omega}\frac{|R^{m-1}g(x,y)|^p}{|x-y|^{mp}}\rho_\varepsilon(|x-y|)\,dx\,dy.
\end{split}
\end{equation}
Combining this with (\ref{eqUpperEstimateLemma}) we arrive at (\ref{eqLemma}).
\end{proof}

The following theorem is an analog of Theorem 2 in \cite{Brezis} (see also \cite{BourgainBrezisMironescu})
for higher order Sobolev spaces.

\begin{theorem}\label{thSobolevSpaces}
Let $\Omega$ be an open set in $\R^n$, $1<p<\infty$ and $m$ be a positive integer.
If  $f\in W^{m-1,p}(\Omega)$ satisfies
\begin{equation}\label{eqMainCondition}
c_f=\liminf\limits_{\varepsilon\to 0}\int_{\Omega}\int_{\Omega}\frac{|R^{m-1}f(x,y)|^p}{|x-y|^{mp}}
\rho_\varepsilon(|x-y|)\,dx\,dy<\infty,
\end{equation}
then $f\in W^{m,p}(\Omega)$.
\end{theorem}

\begin{proof}
Assume that $\Omega'\Subset\Omega$ and $\delta<\dist(\Omega',\partial\Omega)$.
Let $\eta\in C^\infty_0(\R^n)$ be a  nonnegative radial function such that
\[
\int_{\R^n}\eta(x)\,dx=1
\]
and $\supp \eta\subset B(0,1)$. Consider the regularization $f_\delta=f*\eta_\delta$ of $f$, $\eta_\delta(x)=\delta^{-n}\eta(x/\delta)$.
For every $f\in L^1_{\rm loc}(\Omega)$ (extended by zero to $\R^n\setminus\Omega$) the function $f_\delta$ is
smooth in $\Omega$, and if $f$ has a weak derivative $D^\alpha f$ in $\Omega'$, then
$$D^\alpha (f_\delta)=D^\alpha f*\eta_\delta$$ (see e.g. \cite{Mazya}).
Thus, for every $x,y\in \Omega'$ we have
\[
\begin{split}
R^{m-1}&f_\delta(x,y)
=f_\delta(x)-\sum_{|\alpha|\le m-1}D^{\alpha}(f_\delta)(y)\frac{(x-y)^\alpha}{\alpha!}
\\
&=\int_{B(0,\delta)}\bigg(f(x-z)-\sum_{|\alpha|\le m-1}D^{\alpha}f(y-z)\frac{(x-y)^\alpha}{\alpha!}\bigg)\eta_\delta(z)\,dz
\\
&=\int_{B(0,\delta)}R^{m-1}f(x-z,y-z)\eta_\delta(z)\,dz.
\end{split}
\]

By Jensen's inequality, it is easy to see that (\ref{eqMainCondition}) implies
\begin{equation}
\liminf_{\varepsilon\to 0}\int_{\Omega'}\int_{\Omega'}\frac{|R^{m-1}f_\delta(x,y)|^p}{|x-y|^{mp}}\rho_\varepsilon(|x-y|)\,dx\,dy
\le c_f.
\end{equation}

Next by applying (\ref{eqLowerEstimateLemma}) to $g=f_\delta$ we get
\begin{equation}\label{eqSumOfDerivativesEstimate}
\int_{\Omega'}\int_{\partial B(0,1)}\big|\sum_{|\alpha|=m}\frac{D^{\alpha}f_\delta(x)}{\alpha!}e^\alpha\big|^p\,de\,dx
\le c_f.
\end{equation}

Now let $e\in\partial B(0,1)$ and denote by $E$ the vector with the components
$E_\alpha=(1/\alpha!)e_1^{\alpha_1}\cdots e_n^{\alpha_n}$, $|\alpha|=m$. It is easy to see that
$$\Vert v\Vert=\bigg(\int_{\partial B(0,1)} |v \cdot E |^p\,de\bigg)^{1/p},\quad 1\le p<\infty,
$$ is a norm on a linear space of all vectors $v=(v_{\alpha})_{|\alpha|=m}$. Obviously, it is nonnegative and
$\Vert v\Vert=0$ if and only if $v=0$. The triangle inequality follows from the Minkowski inequality.
Since $\Vert\cdot\Vert$ is equivalent to the Euclidean norm, for any vector $\nabla^mg(x)$ we have
\begin{equation}\label{EquivalenceOfNorms}
|\nabla^mg(x)|^p\approx\int_{\partial B(0,1)} |\nabla^mg(x)\cdot E |^p\,de
\end{equation}
and (\ref{eqSumOfDerivativesEstimate}) implies that
\begin{equation}\label{eqGradientOfconvolution}
\int_{\Omega'}|\nabla^mf_\delta(x)|^p\,dx\le c
\end{equation}
for every $\Omega'\Subset\Omega$ and $\delta<\dist(\Omega',\partial\Omega)$,
with the constant $c$ depending only on $c_f$, $n$, $m$ and $p$.

Since $p>1$ the weak compactness and a diagonal argument show that $f\in W^{m,p}(\Omega)$.
\end{proof}

The following statement follows immediately from Theorem \ref{thSobolevSpaces}.

\begin{corollary}
Let $\Omega$ be an open set in $\R^n$, $1\le p<\infty$, $m$ a positive integer.
If $f\in L^p(\Omega)$ satisfies
\begin{equation}\label{eqRecursiveConditions}
\liminf_{\varepsilon\to0}
\int_{\Omega}\int_{\Omega}\frac{|R^{i-1}f(x,y)|^p}{|x-y|^{ip}}\rho_\varepsilon(|x-y|)\,dx\,dy
<\infty,
\end{equation}
for every $i=1,\dots,m$, then $f\in W^{m,p}(\Omega)$.
\end{corollary}

Here the remainders $R^{i-1}f$ are defined recursively starting from $i=1$. Thus, if (\ref{eqRecursiveConditions}) holds for $i=1$ by Theorem \ref{thSobolevSpaces} the weak derivatives $D^\alpha f$, $|\alpha|=1$, exist and $R^{i-1}f$ are defined for $i=2$ and this procedure can be continued recursively.
The recursion may seem somewhat awckward.
A more direct interpretation is possible in terms of Whitney jets, see Section \ref{sWhitneyjets}.

\begin{remark}\label{RmNecessarityOfIntegralConditionP>1}
Note that if $\Omega$ is a $W^{m,p}$-extension domain, i.e there is a bounded linear operator
\[
{\mathcal{E}:W^{m,p}(\Omega)\rightarrow W^{m,p}(\R^n)}
\]
such that $\mathcal{E}f|_{\Omega}=f$ for every $f\in W^{m,p}(\Omega)$, then also the converse statement in Theorem \ref{thSobolevSpaces} is true. Indeed, if $f\in W^{m,p}(\R^n)$, then
\begin{equation}\label{eqLpNormOfRemainderForSobolevFunct}
\int_{\mathbb{R}^n}|R^{m-1}f(x+h,x)|^p\,dx\le c|h|^{mp}\int_{\mathbb{R}^n}|\nabla^mf(x)|^p\,dx
\end{equation}
for every $h\in\R^n$.
This follows easily from Taylor's formula and the fact that smooth functions are dense in the Sobolev space.

Let now $f\in W^{m,p}(\Omega)$ and denote by $\widetilde{f}$ its extension to $\R^n$.
Since
\[
\int_{\mathbb{R}^n}\rho_\varepsilon(|h|)\,dh
=\omega_{n-1}\int_0^\infty\rho_\varepsilon(r)r^{n-1}\,dr
=\omega_{n-1},
\]
where $\omega_{n-1}$ is the $(n-1)$-dimensional surface measure of the
unit ball in $\R^n$,
by (\ref{eqLpNormOfRemainderForSobolevFunct}) we have
\begin{equation}\label{UpperEstimate}
\begin{split}
\int_{\Omega}\int_{\Omega}&\frac{|R^{m-1}f(x,y)|^p}{|x-y|^{mp}}\rho_\varepsilon(|x-y|)\,dy\,dx
\\
&\le\int_{\R^n}\int_{\R^n}\frac{|R^{m-1}\widetilde f(x,y)|^p}{|x-y|^{mp}}\rho_\varepsilon(|x-y|)\,dy\,dx
\\
&\le c\int_{\mathbb{R}^n}|\nabla^m\widetilde{f}(x)|^p\,dx
\\
&\le  c\Vert\widetilde f\Vert_{W^{m,p}(\R^n)}
\le c\Vert f\Vert_{W^{m,p}(\Omega)}.
\end{split}
\end{equation}
\end{remark}

More precisely, the following result is true.

\begin{theorem}\label{thEqualityForSobolevFunction}
Assume that  $\Omega$ is a $W^{m,p}$-extension domain,
let $m$ be a positive integer, $1<p<\infty$ and let $f\in W^{m,p}(\Omega)$.
Then
\begin{equation}\label{eqTheorem1}
\begin{split}
\lim\limits_{\varepsilon\to 0}
\int_{\Omega}\int_{\Omega}&\frac{|R^{m-1}f(x,y)|^p}{|x-y|^{mp}}\rho_\varepsilon(|x-y|)\,dx\,dy
\\
&=\int_{\Omega}\int_{\partial B(0,1)}\Big|\sum_{|\alpha|=m}\frac{D^{\alpha}f(x)}{\alpha! }e^\alpha\Big|^p\,de\,dx.
\end{split}
\end{equation}
\end{theorem}

\begin{proof}
Since $\Omega$ is an extension domain, any function $f\in W^{m,p}(\Omega)$ can be approximated
by functions $f_n\in C^{m+1}_0(\mathbb{R}^n)$ in $W^{m,p}(\Omega)$-norm.
Hence, Lemma \ref{LemmaEqualityForSmoothFunction} implies the validity of (\ref{eqLemma}) for every $f\in W^{m,p}(\Omega)$.
Indeed, by (\ref{EquivalenceOfNorms}) for the right hand side of (\ref{eqLemma}) we have
%If for $g\in W^{m,p}(\Omega)$ we denote by $G$ the vector with the components $D^{\alpha}g(x)/\alpha!$, $|\alpha|=m$, then
\[
\begin{split}
\bigg|\bigg(\int_{\partial B(0,1)}&\Big|\sum_{|\alpha|=m}\frac{D^{\alpha}f(x)}{\alpha! }e^\alpha\Big|^p\,de\bigg)^{1/p}
\\
&\qquad\qquad\qquad
-\bigg(\int_{\partial B(0,1)}\Big|\sum_{|\alpha|=m}\frac{D^{\alpha}f_n(x)}{\alpha! }e^\alpha\Big|^p\,de\bigg)^{1/p}\bigg|
\\
&=|\,\Vert \nabla^mf(x)\Vert-\Vert \nabla^mf_n(x)\Vert\,|
\le\Vert \nabla^m(f-f_n)(x)\Vert
\\
&\le c|\nabla^m(f-f_n)(x)|.
\end{split}
\]
To justify the limit of the left hand side of (\ref{eqLemma}) we can apply (\ref{UpperEstimate}).
Thus, equality (\ref{eqLemma}) is true for any $f\in W^{m,p}(\Omega)$.
%and by (\ref{EquivalenceOfNorms}) we can conclude that the equivalence in (\ref{eqTheorem1}) holds.
\end{proof}

\begin{remark}
If $\Omega=\mathbb{R}^n$ and
\begin{equation}\label{ChoiceOfmollifiers}
\rho_\varepsilon(r)=
\begin{cases}
\dfrac{(n+mp)r^{mp}}{\varepsilon^{n+mp}},&\quad r<\varepsilon,
\\
0,&\quad r\geq\varepsilon,
\end{cases}
\end{equation}
then (\ref{eqMainCondition}) can be written as
\begin{equation}\label{eqMainConditionSpecialCase}
\liminf_{\varepsilon\to 0}\frac{1}{\varepsilon^{mp}}
\int_{\R^n}\mvint_{B(y,\varepsilon)}|f(x)-T^{m-1}_yf(x)|^p\,dx\,dy<\infty.
\end{equation}
Here the integral sign with a bar denotes the integral average.

We point out that  condition (\ref{eqMainConditionSpecialCase}) is closely related to
Calder\'on's characterization of Sobolev spaces in \cite{Calderon1972} (see also \cite{DeVoreSharpley}).
To this end, let $1<p<\infty$ and $m$ be a positive integer.
For $f\in L^p(\Omega)$ we define a maximal function as
\begin{equation}\label{CalderonMaximalFunction}
\mathcal{N}(f,y)=\sup_{\varepsilon>0}\frac{1}{\varepsilon^{m}}
\bigg(\mvint_{B(y,\varepsilon)}|f(x)-P(x,y)|^p\,dx\bigg)^{1/p},
\end{equation}
if there exists a polynomial $P(x,y)$ in $x$, of degree at most $m-1$,
such that the expression on the right hand side of (\ref{CalderonMaximalFunction}) is finite.
If no such polynomial exists, we set $\mathcal{N}(f,y)=\infty$.
Then a function $f\in L^p(\R^n)$ belongs to $W^{m,p}(\R^n)$ if and only if $\mathcal{N}(f,\cdot)\in L^p(\R^n)$.
\end{remark}

Now consider the analog of Theorem \ref{thSobolevSpaces} for $p=1$.
Recall that the space $BV(\Omega)$ is defined as the space of functions in $L^1(\Omega)$, whose weak
derivatives are Radon measures with finite total variation.
Denote by $BV^m(\Omega)$, $m=2,3,\dots$, the set of functions in $L^1(\Omega)$,
whose derivatives of order $m$ are finite Radon measures.

Observe, that by the Riesz representation theorem, a function
$f\in L^1(\Omega)$ belongs to $BV^m(\Omega)$ if and only if there exists $c>0$ such that
\[
\bigg|\int_\Omega fD^\alpha\varphi\,dx\bigg|\le c\Vert\varphi\Vert_{\infty}
\]
for every $\varphi\in C^\infty_0(\Omega)$ and every multi-index $\alpha$ with $|\alpha|=m$.

\begin{theorem}\label{thBVSpaces}
Let $\Omega$ be an open set in $\R^n$ and $m$ be a positive integer. If  $f\in W^{m-1,1}(\Omega)$ satisfies
\begin{equation}\label{eqMainConditionP=1}
\liminf_{\varepsilon\to 0}\int_{\Omega}\int_{\Omega}
\frac{|R^{m-1}f(x,y)|}{|x-y|^{m}}\rho_\varepsilon(|x-y|)\,dx\,dy
<\infty,
\end{equation}
then $f\in BV^m(\Omega)$.
\end{theorem}

\begin{proof}
The proof is the same as for Theorem \ref{thSobolevSpaces}, except for the fact
that (\ref{eqGradientOfconvolution}) implies that all the derivatives of order $m$ of function $f$ are measures with finite total variation.
\end{proof}

\begin{remark}
As observed in \cite{Borghol}, if $\Omega$ is a smooth bounded domain, then $BV^m(\Omega)$
can be characterized as the set of functions
$f\in W^{m-1,1}(\Omega)$ such that $D^\alpha f\in BV(\Omega)$ for every multi-index $\alpha$ with $|\alpha|=m-1$.
The equivalence follows from the Sobolev embedding theorem, see e.g. \cite{Mazya}.

Using this characterization it is not difficult to see that Theorem \ref{thBVSpaces} gives a necessary and sufficient
condition for a function to be in $BV^m(\Omega)$, when $\Omega$ is a bounded smooth domain in $\R^n$.
This can be seen as in Remark \ref{RmNecessarityOfIntegralConditionP>1}.
\end{remark}

\section{A criterion for a function to be a polynomial}

The next result is a higher order version of Theorem 1 in  \cite{Brezis}.

\begin{theorem}\label{thPolynomials}
Let $\Omega$ be a connected open set in $\mathbb{R}^n$, $1\le p<\infty$ and $m$ a positive integer.
Assume that $f\in L^1_{\rm loc}(\Omega)$ has weak derivatives up to order $m-1$ in $\Omega$ and
$$
\liminf_{\varepsilon\to 0}\int_{\Omega}\int_{\Omega}\frac{|R^{m-1}f(x,y)|^p}{|x-y|^{mp}}\rho_\varepsilon(|x-y|)\,dx\,dy=0.
$$
Then $f$ is a polynomial of degree at most $m-1$ a.e. on $\Omega$.
\end{theorem}

\begin{proof}
Suppose first that $p>1$.
Note that if $\Omega$ is a bounded smooth domain in Theorem \ref{thSobolevSpaces}
it is enough to assume that $f\in L^1_{\rm loc}(\Omega)$ has weak derivatives up to order $m-1$ in $\Omega$.
Thus, applying Theorem \ref{thSobolevSpaces}  (see also Theorem \ref{thEqualityForSobolevFunction})
to a ball $B\subset\Omega$ we can conclude that $f\in W^{m,p}(B)$ and $\Vert \nabla^mf\Vert_{L^p(B)}=0$.
This implies that $f$ is a polynomial of degree at most $m-1$ a.e. in $B$
(see e.g. generalized Poincar\'{e} inequality in \cite{Mazya}) and the claim follows from the assumption that
$\Omega$ is connected.

Let then $p=1$.
In this case we apply Theorem \ref{thBVSpaces} to a ball $B\subset\Omega$
and conclude that $f\in BV^m(B)$ and the total variation
$|\nabla^m f|(B)$ of the vector valued measure $\nabla^m f$ equal to zero.
This implies that there is polynomial $P$ of degree at most $m-1$ such that $f=P$ a.e. in $B$
(see e.g. Lemma 12 in \cite{Borghol} for more details).
Again the claim  follows from the assumption that
$\Omega$ is connected.
\end{proof}

The next result is a higher order generalization of Proposition 1 in \cite{Brezis}.

\begin{corollary}\label{cor:polynomial}
Let $\Omega$ be a connected open set in $\mathbb{R}^n$, $1\le p<\infty$ and $m$ a positive integer.
Suppose that a locally integrable function $f$ has weak derivatives up to order $m-1$ in $\Omega$ and
\begin{equation}\label{eq:intfinite}
\int_{\Omega}\int_{\Omega}\frac{|R^{m-1}f(x,y)|^p}{|x-y|^{mp+n}}\,dx\,dy<\infty.
\end{equation}
Then $f$ is a polynomial of degree at most $m-1$ a.e. on $\Omega$.
\end{corollary}

\begin{proof}
By choosing
\[
\rho_\varepsilon(r)=
\begin{cases}
\varepsilon r^{-n+\varepsilon},&\quad r<1,
\\
0,&\quad r\geq 1,
\end{cases}
\]
we have
\begin{equation}\label{eqMainConditionSpecial1}
\begin{split}
&\liminf_{\varepsilon\to0}
\varepsilon\int_\Omega\int_\Omega\frac{|R^{m-1}f(x,y)|^p}{|x-y|^{mp+n-\varepsilon}}\,dx\,dy
\\
&\le\liminf_{\varepsilon\to0}
\varepsilon\int_\Omega\int_\Omega\frac{|R^{m-1}f(x,y)|^p}{|x-y|^{mp+n}}\,dx\,dy
=0.
\end{split}
\end{equation}
The claim follows from Theorem \ref{thPolynomials}.
\end{proof}

Again there is an interpretation of the previous result in terms of Whitney jets.
Indeed, it is possible to state the corollary for Whitney jets without referring to
lower oder derivatives, see Remark \ref{rem:jetremark}.

\section{Whitney jets}\label{sWhitneyjets}

In this section we show that a jet of functions, whose formal Taylor remainder
satisfies (\ref{eqMainConditionSpecialCase}), or its counterpart on a subdomain,
can be identified with the jet of weak derivatives of a Sobolev function.

First we recall terminology related to the Whitney jet theory.
Assume that $\Omega$ is an open set in $\R^n$ and let $m$ be a positive integer.
An $m$-jet $F\in J^m(\Omega)$ is a collection
\[
\{f_\alpha:\,|\alpha|\le m\}
\]
of functions.
The $m$-jets define the formal Taylor polynomials in $x\in\mathbb{R}^n$ (centered at $y\in \Omega$)
$$
T_y^{k}F(x)=\sum_{|\alpha|\le k}f_{\alpha}(y)\frac{(x-y)^\alpha}{\alpha!},
$$
with $k\le m$, and
$$
T_{y,j}^{k-|j|}F(x)=\sum_{|j+\alpha|\le k}f_{j+\alpha}(y)\frac{(x-y)^\alpha}{\alpha!},
$$
with $|j|\le k\le m$.
The formal Taylor remainders are defined to be
\[
R^{k}F(x,y)=f_0(x)-T_y^{k}F(x)
\]
and
\[
R^{k-|j|}_{j}F(x,y)=f_j(x)-T_{y,j}^{k-|j|}F(x),
\]
where $x,y\in K$ and $|j|\le k\le m$.

Let $Q$ be a fixed cube in $\R^n$, an $(m-1)$-jet $F\in J^{m-1}(Q)$,
\[
F=\{f_j:|j|\le m-1\},
\]
where $f_j\in L^p(Q)$, is said to be an $(m-1)$-jet in $Q$ with variable Lipchitz coefficients,
denoted by $F\in \vlc(m,p,Q)$, if the pointwise inequality
\begin{equation}\label{VLCcondition}
|R^{m-1}F(x,y)|\le|x-y|^m(a_Q(x)+a_Q(y)),
\end{equation}
where $x,y\in Q$, holds for some function $a_Q=a_Q(F)\in L^p(Q)$.

The $(m-1)$-jet spaces $\vlc(m,p,Q)$ have been studied in \cite{BojarskiProcSteklovInst} and it has been shown that $\vlc(m,p,Q)$
regarded as a Banach space and equipped with the norm
$$\Vert F\Vert=\max\{\Vert f_j\Vert_{L^p(Q)}:|j|\le m-1\}+\inf\Vert a_Q\Vert_{L^p(Q)},$$
can be identified with the classical Sobolev space $W^{m,p}(Q)$.
The fact that inequality (\ref{VLCcondition}) holds for a function $f\in W^{m,p}(Q)$ has been proved before in \cite{BojarskiHajlaszStrelecki}.

Let us consider an $(m-1)$-jet $F\in J^{m-1}(\Omega)$ of locally integrable functions on an open set $\Omega$ in $\R^n$
with the property that  its formal Taylor remainder of order $m-1$ satisfies the condition
\begin{equation}\label{eqMainConditionTh4}
a_F=\lim_{\varepsilon\to 0}\frac{1}{\varepsilon^{n+mp}}\iint_{\{(x,y)\in\Omega:|x-y|<\varepsilon\}}
|R^{m-1}F(x,y)|^p\,dx\,dy<\infty.
\end{equation}
This is a special case of condition (\ref{eqMainCondition}) with mollifiers (\ref{ChoiceOfmollifiers}).

Note that pointwise estimate (\ref{VLCcondition}) for $x,y\in\Omega$ implies (\ref{eqMainConditionTh4}).
Indeed,
\[
\begin{split}
\int_\Omega\int_\Omega& |R^{m-1}F(x,y)|^p\chi_{\{|x-y|<\varepsilon\}}\,dx\,dy
\\
&\le\int_\Omega\int_\Omega |x-y|^{mp}\chi_{\{|x-y|<\varepsilon\}}(a_{\Omega}(x)+a_{\Omega}(y))^p\,dx\,dy
\\
&\le c\int_\Omega\int_\Omega |x-y|^{mp}\chi_{\{|x-y|<\varepsilon\}}a^p_{\Omega}(x)\,dx\,dy
\end{split}
\]
and we have
\[
\begin{split}
\frac{1}{\varepsilon^{n+mp}}
&\iint_{\{(x,y)\in\Omega:|x-y|<\varepsilon\}} |R^{m-1}F(x,y)|^p\,dx\,dy
\\
&\le c\int_\Omega a^p_{\Omega}(x) \int_\Omega\frac{|x-y|^{mp}}{\varepsilon^{n+mp}}\chi_{\{|x-y|<\varepsilon\}}\,dy\,dx
\\
&\le c\frac{\omega_{n-1}}{n}\Vert a_\Omega\Vert^p_{L^p(\Omega)}.
\end{split}
\]

At the same time (\ref{eqMainConditionTh4}) is a sufficient condition for a jet to be identified with a Sobolev function.
More precisely, the following theorem holds true.

\begin{theorem}\label{thWhitneyjet}
Let $\Omega$ be an open set in $\mathbb{R}^n$, $1\le p<\infty$ and $m$ be a positive integer.
Assume that an $(m-1)$-jet
\[
F=\{f_\alpha:\,|\alpha|\le m-1\},
\]
where $f_\alpha\in L^p(\Omega)$, satisfies condition (\ref{eqMainConditionTh4}).
Then for every $\Omega'\Subset\Omega$ there is a function $f\in W^{m,p}(\Omega')$ such that
\[
f_\alpha|_{\Omega'}=D^{\alpha}f|_{\Omega'},\,\,\,\,|\alpha|\le m-1,
\]
and
\[
\|\nabla^m f\|_{L^p(\Omega')}\le c\, a_F.
\]
\end{theorem}

\begin{remark}\label{rem:jetremark}
Corollary \ref{cor:polynomial} reads for Whiney jets as follows.
Let $\Omega$ be a connected open set in $\mathbb{R}^n$.
Suppose that
\[
\int_{\Omega}\int_{\Omega}\frac{|R^{m-1}F(x,y)|^p}{|x-y|^{mp+n}}\,dx\,dy<\infty.
\]
Then $f$ is a polynomial of degree at most $m-1$ a.e. on $\Omega$.
\end{remark}

To prove the theorem we use the sketch of the proof of Theorem 9.1 from \cite{BojarskiProcSteklovInst}.
First we show that the next statement is true.

\begin{lemma}\label{LemmaConditionsForJet}
Let $\Omega$ be an open set in $\mathbb{R}^n$, $1\le p<\infty$ and $m$ be a positive integer.
Suppose that $(m-1)$-jet $F$ satisfies condition (\ref{eqMainConditionTh4}).
Then for every $\Omega'\Subset\Omega$ we have
\begin{equation}\label{eqMainConditionForJet}
\lim_{\varepsilon\to 0}\frac{1}{\varepsilon^{(m-|j|)p}}\int_{\Omega'}
\mvint_{B(x,\varepsilon)}|R_j^{m-1-|j|}F(x,y)|^p\,dy\,dx\le c\,a_F,
\end{equation}
whenever $|j|\le m-1$.
\end{lemma}

\begin{proof}
Let $\Omega'\Subset\Omega$ and $0<\varepsilon<\dist(\Omega',\partial\Omega)$. Fix $x,y\in\Omega'$, $|x-y|<\varepsilon$. Using Taylor algebra arguments, we have
\begin{equation}\label{RemainderViaPolynomialDerivative}
\begin{split}
R^{m-1-|j|}_j F(x,y)
&=D^j_z[R^{m-1}F(z,x)-R^{m-1}F(z,y)]_{z=x}
\\
&= D_z^jP(z;x,y)_{z=x},
\end{split}
\end{equation}
where $P(z;x,y)$ is a polynomial in $z$ of order at most $m-1$.

Since  $\varepsilon<\dist(\Omega',\partial\Omega)$ the set
$S=B(x,\varepsilon)\cap B(y,\varepsilon)\subset\Omega.$
It is easy to see that
$$c|S|\geq |B(x,\varepsilon)|=|B(y,\varepsilon)|$$
for some constant $c$ which is independent of $\varepsilon$.
By Markov's inequality \cite{Campanato} applied to the subset $S$ of the ball $B(x,\varepsilon)$, we obtain
\begin{equation*}%\label
|D_z^jP(z;x,y)|_{z=x}\le\frac{c(n)}{\varepsilon^{|j|}}\bigg(\mvint_{S}|P(x';x,y)|^p\,dx'\bigg)^{1/p}.
\end{equation*}

Thus, from (\ref{RemainderViaPolynomialDerivative}), we have
\begin{equation}\label{ProofLemma}
\begin{split}
|R^{m-1-|j|}_j F(x,y)|
&\le\frac{c(n)}{\varepsilon^{|j|}}\bigg[\bigg(
\mvint_{S}|R^{m-1}F(x',x)|^p\,dx'\bigg)^{1/p}
\\
&\qquad\qquad+\bigg(\mvint_{S}|R^{m-1}F(x',y)|^p\,dx'\bigg)^{1/p}\bigg]
\\
&=\frac{c(n)}{\varepsilon^{|i|}}\big[I(x)^{1/p}+I(y)^{1/p}\big].
\end{split}
\end{equation}

It is easy to see that
\[
\begin{split}
\int_{\Omega'}\mvint_{B(x,\varepsilon)}I(x)\,dy\,dx
&\le c\int_{\Omega'}\mvint_{B(x,\varepsilon)}\mvint_{B(x,\varepsilon)}
|R^{m-1}F(x',x)|^p\,dx'\,dy\,dx
\\
&\le c\iint_{\{(x,x')\in\Omega:|x-x'|<\varepsilon\}}|R^{m-1}F(x',x)|^p\,dx'\,dx.
\end{split}
\]
On the other hand, we obtain a similar estimate for
\[
\int_{\Omega'}\mvint_{B(x,\varepsilon)}I(y)\,dy\,dx,
\]
which together with (\ref{ProofLemma}) proves the claim.
\end{proof}

\begin{remark} Since  (\ref{eqMainConditionForJet}) implies
\begin{equation*}\label{eqMainConditionForJet1}
\frac{1}{\varepsilon^{n+(m-|j|)p}}\iint_{\{(x,y)\in\Omega':|x-y|<\varepsilon\}}
|R_j^{m-1-|j|}F(x,y)|^p\,dx\,dy\le c\,a_F,
\end{equation*}
as $\varepsilon\to 0$, $|j|\le m-1$,
using the terminology of Jonsson and Wallin (see e.g. \cite{JonssonWallin}) we can formulate Lemma \ref{LemmaConditionsForJet} in the following way:
If $(m-1)$-jet $F$ satisfies condition (\ref{eqMainConditionTh4}) then for every $\Omega'\Subset\Omega$
we have $F\in{\rm Lip }(m,p,\infty,\Omega')$.
%$$=
%\frac{(2c(n))^p}{\varepsilon^{|i|p}}\int\limits_{\Omega'}\int\limits_{\Omega'}\chi_{\{|x-y|<\varepsilon\}}\mathop{\int\limits_{S}\mspace{-19
%mu}{-}}|R^{m-1}F(x',x)|^p\,dx'\,dxdy\le \frac{c}{\varepsilon^{|i|p}}\underset{(x,y)\in\Omega:|x-y|<\varepsilon}{\int\int}|R^{m-1}F(x',x)|^p\,dx'\,dx.
%$$
\end{remark}

\begin{proof}[Proof of Theorem \ref{thWhitneyjet}]
%{\it Proof of Theorem \ref{thWhitneyjet}.}
Let $\Omega'$ be an open set such that $\Omega'\Subset\Omega$. Decompose $\Omega'$ into dyadic cubes.
More precisely, let $\mathcal{M}_k$ denote a net with mesh $2^{-k}$ in $\mathbb{R}^n$ i.e. $\mathcal{M}_k$ is a division of $\mathbb{R}^n$ into equally large closed cubes with side lenghts $2^{-k}$, obtained by slicing $\mathbb{R}^n$ with hyperplanes orthogonal to the coordinate axis.

Set $$\Omega'_k=\big\{ x\in \Omega':\,{\rm dist}(x,\Omega'^c)>\sqrt{n}2^{-k+1}\big\}$$
and
$$\Omega'_l=\big\{ x\in \Omega':\,\sqrt{n}2^{-l+1}<{\rm dist}(x,\Omega'^c)\le\sqrt{n}2^{-l+2}\big\},$$
for $l=k+1,\dots$
Define
$$\mathcal{F}^0_k=\bigcup\limits_{l=k}^{\infty}\{Q\in \mathcal{M}_{l}:\,Q\cap\Omega'_l\neq\emptyset\}$$
and denote by $\mathcal{F}_k$ the collection of maximal cubes of $\mathcal{F}^0_k$
(see e.g. \cite{Stein} for details on the Whitney decomposition).

Then, for each $k$ the collection of cubes $\mathcal{F}_k=\{Q_i^k\}_{i\in I_k}$ satisfies the conditions:
\begin{itemize}
\item[(i)] $\Omega'=\bigcup\limits_{i\in I_k}Q_i^k$;
\item[(ii)]${\rm diam} Q_i^k\le 2^{-k}\sqrt{n}$ for all $i\in I_k$, moreover,
\[
{\rm diam} Q_i^k\le {\rm dist}(Q_i^k,\Omega'^c)\le 4{\rm diam}Q_i^k,
\]
if ${\rm diam} Q_i^k\le 2^{-k-1}\sqrt{n}$;
\item[(iii)] the interiors of the cubes $Q_i^k$ are pairwise disjoint.
\end{itemize}

It is easy to see that theses properties imply the next statements as well:
\begin{itemize}
\item[(iv)] $\sum\limits_{i\in I_k}\chi_{\frac{9}{8}Q_i^k}\le c$, where constant $c$ depends only on $n$;
\item[(v)] if $Q,\,K\in\mathcal{F}_k$ and $\frac{9}{8}Q\cap\frac{9}{8}K\neq\emptyset$,
then ${\rm diam} Q\le {\rm diam} K\le 4{\rm diam}Q$.
\end{itemize}

Let now $\{\varphi_i^k:i\in I_k\}$ be a smooth partition of unity subordinated to the decomposition $\mathcal{F}_k$.
It's standard properties (see \cite{Stein}) are:
\begin{itemize}
\item[(i)]  ${\rm supp}\varphi_i^k\subset \mathcal{Q}^i_k=\frac{9}{8}Q_i^k$,
\item[(ii)] $\sum\limits_{i\in I_k}\varphi_i^k\equiv 1$ on $\Omega'$ and
\item[(iii)] $|D^\alpha\varphi_i^k|\le c_\alpha ({\rm diam} Q_i^k)^{-\alpha}$ in $\mathbb{R}^n$.
\end{itemize}

Define the approximating family
\begin{equation}
w^k(x)=\sum\limits_{i\in I_k}\varphi_i^k(x)T_{\mathcal{Q}_i^k}^{m-1}F(x),
\end{equation}
where $x\in \Omega'$ and
\[
T_Q^{m-1}F(x)=\mvint_QT^{m-1}_yF(x)\,dy.
\]
Let $x$ be in $\Omega'$, denote by $B^k_x$ the ball $B(x,r_k)$,
\[
r_k=\sqrt{n}\max\{{\rm diam}\,\mathcal{Q}_i^k:\,x\in\mathcal{Q}_i^k\}.
\]
By property (v) of cubes $Q_i^k$ the radius $r_k\approx{\rm diam}\,Q_i^k$ whenever $x\in{Q}_i^k$ and, evidently, $\mathcal{Q}_i^k\subset B_x^k$ for every $i$ such that $x\in \mathcal{Q}_i^k$. Consider
\[
\begin{split}
|f_0(x)-w^k(x)|
&=\bigg|\sum\limits_{i\in I_k}\varphi_i^k(x)
\mvint_{\mathcal{Q}_i^k}[f_0(x)-T^{m-1}_yF(x)]\,dy\bigg|
\\
&\le\sum\limits_{i\in I_k}\varphi_i^k(x)\mvint_{\mathcal{Q}_i^k}
|R^{m-1}F(x,y)|\,dy\,\chi_{\mathcal{Q}_i^k}(x)
\\
&\le c\sum\limits_{i\in I_k}\varphi_i^k(x)\mvint_{B_x^k}|R^{m-1}F(x,y)|\,dy
\\
&=c\mvint_{B_x^k}|R^{m-1}F(x,y)|\,dy.
\end{split}
\]
Thus,
$$
\Vert f_0-w^k\Vert_{L^p(\Omega')}\le c\bigg(\int_{\Omega'}\mvint_{B_x^k}|R^{m-1}F(x,y)|^p\,dy\,dx\bigg)^{1/p}\le c2^{-km},
$$
as $k\to\infty$,
which shows that the sequence of functions $w^k$ converges in $L^p(\Omega')$ to $f_0$.

Now let $|\alpha|=m$ and consider
\[
\begin{split}
D^\alpha w^k(x)
&=\sum\limits_{i}D^\alpha\big(\varphi_i^k(x)T_{\mathcal{Q}_i^k}^{m-1}F(x)\big)
\\
&=\sum\limits_{i}\sum\limits_{\beta\le \alpha}
\frac{\alpha!}{\beta!(\alpha-\beta)!}
D^\beta\varphi_i^k(x)D^{\alpha-\beta}(T_{\mathcal{Q}_i^k}^{m-1}F)(x)
\\
&=\underset{
|\gamma|\le m-1}{\sum\limits_{\beta+\gamma=\alpha}}\frac{\alpha!}{\beta!(\alpha-\beta)!}
S_{\beta,\gamma}(x),
\end{split}
\]
where
$$
S_{\beta,\gamma}(x)=\sum\limits_{i}D^\beta\varphi_i^k(x)T_{\mathcal{Q}_i^k}^{m-1-|\gamma|}D^\gamma F(x).
$$
Since $\sum_iD^\beta\varphi_i(x)\equiv 0$ if $|\beta|\geq 1$,
\[
\begin{split}
S_{\beta,\gamma}(x)
&=S'_{\beta,\gamma}(x)\equiv\sum\limits_{i}D^\beta\varphi_i^k(x)\big[f_\gamma(x)-T_{\mathcal{Q}_i^k}^{m-1-|\gamma|}D^\gamma F(x)\big]
\\
&=\sum\limits_{i}
D^\beta\varphi_i^k(x)\mvint_{\mathcal{Q}_i^k}R_\gamma^{m-1-|\gamma|}F(x,y)\,dy
\end{split}
\]
and
$$
|S'_{\beta,\gamma}(x)|\le\sum\limits_{i}
|D^\beta\varphi_i^k(x)|\mvint_{B_x^k}|R_\gamma^{m-1-|\gamma|}F(x,y)|\,dy\,\chi_{\mathcal{Q}_i^k}(x).
$$
By smoothness properties of the partition of unity
\[
|D^\beta\varphi_i^k(x)|\le c({\rm diam}\,Q_i^k)^{-|\beta|}\le c 2^{k|\beta|},
\]
for every $x\in \mathcal{Q}_i^k$. Hence, the $L^p$-norm of $S'_{\beta,\gamma}$ can be estimated in the following way
\[
\begin{split}
\Vert S'_{\beta,\gamma}\Vert_{L^p(\Omega')}^p
&=\int_{\Omega'}\bigg(\sum\limits_{i}
|D^\beta\varphi_i^k(x)|\mvint_{B_x^k}
|R_\gamma^{m-1-|\gamma|}F(x,y)|\,dy\,\chi_{\mathcal{Q}_i^k}(x)\,\bigg)^p\,dx
\\
&\le c2^{kp|\beta|}\int_{\Omega'}\bigg(
\mvint_{B_x^k}|R_\gamma^{m-1-|\gamma|}F(x,y)|\,dy\,\bigg)^p\,dx
\\
&\le c2^{kp|\beta|}\int_{\Omega'}
\mvint_{B_x^k}|R_\gamma^{m-1-|\gamma|}F(x,y)|^p\,dy\,dx.
\end{split}
\]
Thus, by Lemma \ref{LemmaConditionsForJet} for $|\alpha|=m$ and large $k$ we have
\begin{equation}\label{eqMthDerivative}
\Vert D^\alpha w^k\Vert_{L^p(\Omega')}\le c a_F 2^{(|\beta|+|\gamma|-m)k}\le c a_F
\end{equation}
(one could set $\varepsilon=2^{-k}$ in (\ref{eqMainConditionForJet})).

We show that for $|\alpha|\le m-1$ the derivatives $D^\alpha w^k$ converge in $L^p(\Omega')$ to
the elements $f_\alpha$ of the jet $F$.
A computation shows that
\[
\begin{split}
f_\alpha(x)-D^\alpha w^k(x)
&=\sum\limits_{i}\varphi_i^k(x)f_\alpha(x)-\sum\limits_{i}D^\alpha\big(\varphi_i^k(x)T_{\mathcal{Q}_i^k}^{m-1}F(x)\big)
\\
&=\sum\limits_{i}\varphi_i^k(x)\big(f_\alpha(x)-T_{\mathcal{Q}_i^k}^{m-1-|\alpha|}D^\alpha F(x)\big)
\\
&\qquad+\underset{|\beta|>0}{\sum\limits_{\beta+\gamma=\alpha}}
\frac{\alpha!}{\beta!(\alpha-\beta)!}\sum\limits_{i}
D^\beta\varphi_i^k(x)T_{\mathcal{Q}_i^k}^{m-1-|\gamma|}D^\gamma F(x)
\\
&=\sum\limits_{\beta+\gamma=\alpha}
\frac{\alpha!}{\beta!(\alpha-\beta)!}S'_{\beta,\gamma}(x).
\end{split}
\]
Using the above estimate for $\Vert S'_{\beta,\gamma}\Vert_{L^p(\Omega')}^p$, we arrive at
\begin{equation}\label{ConvergenceOfJetToDerivatives}
\Vert f_\alpha- D^\alpha w^k\Vert_{L^p(\Omega')}\le c 2^{(|\beta|+|\gamma|-m)k}= c 2^{(|\alpha|-m)k}\to 0,
\end{equation}
as $k\to\infty$.

Thus, we have shown that the sequence of smooth functions $w^k$ converges to function $f_0$ in $L^p(\Omega')$ and for any $\alpha$, $|\alpha|\le m$, the sequence of derivatives $D^\alpha w^k$ is bounded in $L^p(\Omega')$. Since $p>1$, by weak compactness argument  it follows that $f_0\in W^{m,p}(\Omega')$.

Let $f\equiv f_0$, then  by \eqref{ConvergenceOfJetToDerivatives} $D^{\alpha}f|_{\Omega'}=f_\alpha|_{\Omega'}$, $|\alpha|\le m-1$, and by \eqref{eqMthDerivative} we have estimate $\|\nabla f^m\|_{L^p(\Omega')}\le c a_F$.
\end{proof}
%$$\Vert f-w^k\Vert_{W^{m-1,p}(\Omega)}\to 0\,\,\,\, \text{as}\,\,\,k\to\infty.$$

\begin{remark}
We can also consider the meaning of the condition (\ref{eqMainConditionTh4}) for the Whitney
jet on an $s$-set in $\R^n$.
Recall that a subset $S\subset\mathbb{R}^n$ is called an $s$-set (or Ahlfors $s$-regular)
if there are constants $c_1,\,c_2>0$ such that for every ball $B=B(x,r)$ with
center at $S$ and $r\le {\rm diam}S$
\[
c_1r^s\le H^s(B(x,r)\cap S)\le c_2r^s,
\]
where $H^s$ denote the $s$-dimensional Hausdorff measure on $\mathbb{R}^n$.

Let $S$ be an $s$-set with $n-1< s\leq n$. Consider a $(m-1)$-jet
\[
F=\{f_\alpha:\,|\alpha|\le m-1\},
\]
$f_\alpha\in L^p(S,H^s)$,
such that it's formal Taylor remainder of order $m-1$ satisfies the condition
\begin{equation}\label{eqMainConditionSec5}
\frac{1}{\varepsilon^{s+mp}}\iint_{\{(x,y)\in S:|x-y|<\varepsilon\}}|R^{m-1}F(x,y)|^p\,dH^s\,dH^s
\le a_F
\end{equation}
as $\varepsilon\to 0$.
Then the functions $f_\alpha$ are uniquely determined by $f_0$. The proof uses slightly modified Lemma \ref{LemmaConditionsForJet} and arguments as in Theorem 2 in \cite{JonssonWallin}, p.126. Thus, if $f_0$ is identically zero then  $f_\alpha$ are zero as well and we can identify the jet $F$ with the single function $f_0$.

In particular, if $\Omega$ is a $W^{m,p}$-extension domain then $\Omega$ is an $n$-set (see \cite{HajlaszKoskelaTuominen}) and the statement '' a jet $F$ satisfies (\ref{eqMainConditionSec5}) on $\Omega$'' is equivalent to the fact that $f_0\in W^{m,p}(\Omega)$. If $S$ is an arbitrary $s$-set with $n-1< s\leq n$ then condition (\ref{eqMainConditionSec5}) characterizes the Lipschitz space ${\rm Lip }(m,p,\infty,S)$ introduced by A. Jonsson and H. Wallin \cite{JonssonWallin}.

The first order spaces ${\rm Lip }(1,p,\infty,S)$ and their relations to Sobolev spaces on metric measure spaces were studied in \cite{YangLin}. For instance, from the results in \cite{YangLin} it follows that if $S$ supports the $q$-Poncar\'{e} inequality, $1\le q<\infty$, then ${\rm Lip }(1,p,\infty,S)$ coincides with the Haj\l asz-Sobolev space $W^{1,p}(S)$ for $p>q$.
\end{remark}

\bigskip
\noindent Addresses:

\smallskip
\noindent B.B.: Institute of Mathematics, Polish Academy of Sciences,
00-956, Warsaw, Poland.\\
\noindent
E-mail: {\tt b.bojarski@impan.gov.pl}\\

\noindent L.I.: Department of Mathematics, P.O. Box 11100,
FI-00076 Aalto University, Finland. \\
\noindent
E-mail: {\tt liza.ihnatsyeva@tkk.fi}\\

\noindent J.K.: Department of Mathematics, P.O. Box 11100,
FI-00076 Aalto University, Finland. \\
\noindent
E-mail: {\tt juha.kinnunen@tkk.fi}\\

\end{document}